\documentclass[12pt]{article} 
\usepackage[T1]{fontenc}
\usepackage[latin1]{inputenc}
\usepackage{amssymb}
\usepackage{latexsym}
\usepackage{amsmath}

\usepackage[left=3cm,top=4cm,bottom=3cm,right=3cm]{geometry}

\usepackage{graphicx}
\input xy
\input diagxy
\xyoption{all}

\newcommand{\vf}[0]{{\varphi}}

\newcommand{\memof}[0]{\, \epsilon \,}
\newcommand{\subseteqof}[0]{\, \dot{\subseteq} \,}

\newcommand{\mono}[0]{\to/ >->/}
\newcommand{\pto}[0]{\rightharpoondown}

\newcommand{\idmap}[0]{{\sf id}}
\newcommand{\domain}[0]{{\sf dom}}
\newcommand{\codomain}[0]{{\sf cod}}
\newcommand{\fst}[0]{{\sf fst}}
\newcommand{\snd}[0]{{\sf snd}}
\newcommand{\comp}[0]{{\sf comp}}
\newcommand{\bool}[0]{{\bf 2}}
\newcommand{\terminal}[0]{{\bf 1}}
\newcommand{\initial}[0]{{\bf 0}}

\newcommand{\piev}[0]{{\rm ev}}

\newcommand{\qed}[0]{\Box}

\newcommand{\omitthis}[1]{}

\newcommand{\longtext}[1]{}
\newcommand{\shorttext}[1]{}

\newtheorem{thm}{Theorem}[section]
\newtheorem{cor}[thm]{Corollary}
\newtheorem{prop}[thm]{Proposition}
\newtheorem{lemma}[thm]{Lemma}

\newtheorem{defin}[thm]{Definition}

\newtheorem{remarks}[thm]{Remarks}

\title{Constructivist and Structuralist Foundations: \\
Bishop's and  Lawvere's Theories  of Sets }
\author{
Erik Palmgren\footnote{Research supported by
 Swedish Collegium for Advanced Study, Uppsala,
Swedish Research Council (VR) and Institut Mittag-Leffler, 
Djursholm, Sweden.
}}

\date{October 23, 2009. Revised: August 11, 2010 and July 14, 2011}

\begin{document}

\maketitle

\begin{abstract} Bishop's informal set theory is briefly discussed and
compared to Lawvere's Elementary Theory of the Category of Sets (ETCS).
We then present a constructive and predicative version of ETCS, whose
standard model is based on the constructive type theory of
Martin-L{\"o}f. The theory, CETCS, provides a structuralist foundation
for constructive mathematics in the style of Bishop.


\end{abstract}

\section{Introduction}

Errett Bishop's book {\em Foundations of Constructive Analysis} 
from 1967 contains a chapter on set theory. This set theory, apart
from being informal, is quite unlike any of the theories
of Zermelo--Fraenkel or G{\"o}del--Bernays, which are derived from
the iterative concept of set.
\begin{quote}
``A set is not an entity which has an ideal existence: a set exists only
when it has been defined. To define a set we prescribe, at least
implicitly, what we (the constructing intelligence) must do in order
to construct an element of the set, and what we must do to show that
two elements are equal.'' (Bishop 1967, p.\ 2)
\end{quote}
We find a similar explanation of what a set is also in the type theory
of Martin-L{\"o}f (1984). Both explanations are aligned to Cantor's
early explanation of sets from 1882 in the respect that they mention
conditions for equality of elements explicitly. See Tait (2000) for a
discussion.  Bishop (1967, p.\ 74) emphasizes
that two elements may not be compared unless they belong to some
common set. This indicates a type-theoretic attitude to the
foundations.  Bishop's version of set theory has, despite its
constructiveness, a more abstract character than e.g.\ ZF set theory
in that it does not concern coding issues for basic mathematical
objects. It defines a subset of a set $X$ to be a pair $(A,i_A)$ where
$i_A:A \to X$ is function so that $a = b$ if, and only if, $i(a)=
i(b)$.  An element $x\in X $ is a member of the subset if $x=i_A(a)$
for some $a \in A$. That the subset $(A,i_A)$ is included in another
subset $(B,i_B)$ of $X$ is defined by requirement that there is a
function $f:A \to B$ so that $i_A= i_B \circ f$, i.e.\ that the
diagram
\begin{equation} \label{bishincl}
\bfig
\Vtriangle[A`B`X;f`i_A`i_B]
\efig
\end{equation}
commutes. The subsets are equal in case $f$ is a bijection.
Unions and intersection are only defined when the involved sets are subsets of
the same underlying set. These and other features of Bishop's set theory 
are remarkably reminiscent of Lawvere's {\em Elementary Theory of the
Category of Sets} (ETCS) introduced in 1964. ETCS is obtained by
singling out category-theoretic universal properties of various set
construction in such a way that they become invariant under
isomorphism; see  (McLarty 2004) and the introduction (McLarty 2005)
 to (Lawvere 2005), the
full version of the 1964 paper. This invariance is of course fundamental 
for a {\em structuralist foundation}. ETCS is an elementary theory 
in the sense that it uses classical first order logic as a basis, 
and make no special assumption on existence of second order or higher 
order objects. The theory is equivalent to the axioms of a well-pointed
 topos with the axiom of choice (McLarty 2004, MacLane 1998). 
It should be emphasized that ETCS was introduced to give an immediate
axiomatization of sets, while the Lawvere-Tierney elementary theory 
of a topos was intended to give axioms for sheaves of sets over an arbitrary
topological space.

Bishop (1970a, 1970b) considered various versions of G{\"o}del's system T
 as a possible foundation for his set
theory. At the basis of the interpretation is a system of computable
functions and functionals, which in effect are the core operations of certain
modern programming languages.
Full-fledged systems suitable for the formalization of
constructive mathematics in the style of Bishop emerged later with the
constructive type theory of  Martin-L{\"o}f (1975) and the constructive set
theories CST (Myhill 1975) and CZF (Aczel 1978). Of these, the type-theoretic 
system is the more fundamental from a constructive semantical point of view, since it
describes explicitly how the computation of functions are carried out.
Indeed, the mentioned set-theoretic system, CZF, can be justified on the
grounds of Martin-L{\"o}f's type theory (MLTT) as shown by Aczel (1978)
by a model construction. In MLTT the explanation of when elements of
a set (type) are equal halts at the level of definitional equality. There are
no quotient constructions, so it is customary to consider  a  type together
with an equivalence relation, as a set-like object, a so-called {\em setoid}.
This gives two possible conceptions of constructive sets based on the 
formal theories CZF and MLTT, namely iterative sets 
(sets as trees) and setoids respectively.

In this paper we present a constructive version of ETCS, called CETCS,
which is obtained abstracting on category-theoretic properties of CZF
sets and of setoids in a universe in MLTT.  A first requirement on
CETCS is of course that we use intuitionistic first order logic
instead of the customary classical logic.  CETCS has however the
property that by adding the law of excluded middle and the axiom of
choice (AC), we get a theory equivalent to ETCS.  Furthermore the
theories of Aczel--Myhill and Martin-L{\"o}f are (generalized)
predicative, so that power set principles are
not valid. Thus a constructive ETCS cannot be obtained by adding
axioms to the elementary theory of toposes. In Moerdijk and Palmgren
(2000, 2002) a notion of predicative topos was introduced taking the setoids
of MLTT with a hierarchy of
universes as a standard model.  Other variants of predicative toposes
have been introduced and studied (van den Berg 2005); see also 
Maietti (2005) and Awodey and Warren (2005). A drawback of
the category of setoids, as opposed categories of sets, is that there
is no canonical choice of pullbacks (Sect.\ 6, Hofmann 1994).  This
makes the formulation of some axioms a bit less concise, but also more
general. 

We emphasize that ETCS does not deal with the
set-class distinction or replacement axioms. ETCS with replacement
has however been considered (Osius 1974, McLarty 2004). A constructive
treatment of the set-class distinction was given by Joyal and Moerdijk (1995)
by the introduction of notion of a small map. Predicatively acceptable versions
of this were developed in (Moerdijk and Palmgren 2002) and 
(Moerdijk and van den Berg 2008). It seems rather straightforward
to extend CETCS to include axioms for small maps along those lines. 
Another possible extension of CETCS is to add inductively defined
subsets. We leave these investigations for another occasion.
A feature of CETCS is that it introduces a constructive version of well-pointedness. Shulman (2010)
gives a definition of this notion which works for weaker categories.

An outline of the paper is a follows: In Section 2 a standard
first-order logic definition of categories is given. We present in
Section 3 some notation regarding relations and subobjects for
categories where products are not supposed to be chosen.  The axioms
of ETCS and CETCS are presented in parallel and compared in Section 4.
In Section 5 some elementary set-theoretic consequence are drawn from
CETCS, which indicates its usefulness for Bishop style constructive
mathematics. It is shown that CETCS together with the axiom of choice
and classical logic gives the original ETCS. The relation of CETCS to
standard category theory notions is given in Section 6 and Section
7. This can part can be skipped by the reader that is not 
particularly interested
in categorical logic.  Section 7 contains a technical contribution 
which shows how a ``functor free'' formulation of locally
cartesian closed categories (LCCCs) can
be employed in categorical logic.

\subsubsection*{Acknowledgment}
The main results of this article were obtained while the author was a
fellow of the Swedish Collegium for Advanced Study, January through
June 2009.  Many thanks go to the Collegium and its principal
Professor Bj{\"o}rn Wittrock for the opportunity to work in this
stimulating research environment. The author is grateful to 
Institut Mittag-Leffler for support, and to the anonymous referee 
for helpful remarks.

% ...

\section{Elementary Categories}

We shall take care to formulate all the axioms so that they may be easily cast in many sorted first-order (intuitionistic) logic.
Following the notation of Mac Lane (1998),
a category ${\cal C}$ is specified by an algebraic signature consisting of three
collections ${\cal C}_0,{\cal C}_1,{\cal C}_2$ (for objects, mappings
(or arrows),
composable mappings) and six functions 
$\idmap: {\cal C}_0 \to {\cal C}_1$, 
$\domain, \codomain: {\cal C}_1 \to {\cal C}_0$,
$\comp: {\cal C}_2 \to {\cal C}_1$, $\fst,\snd: {\cal C}_2 \to {\cal C}_1$.
The intention is that $\domain$ gives the domain of the mapping while $\codomain$ gives its codomain. 
The collection ${\cal C}_2$ is supposed to consist of composable mappings
$$\cdot \to^f \cdot \to^g \cdot$$
and $\fst$ gives the first of these mappings while  $\snd$ gives the
second mapping. Then $\comp$ is the composition operation.
The axioms for a category are then briefly as follows, where variables ranges
are $x \in {\cal C}_0$, $f,g,h,k,\ell \in {\cal C}_1$, $u,v \in {\cal C}_2$:
(K1) $\domain(\idmap_x) = x$,  (K2) $\codomain(\idmap_x) = x$,
 (K3) $\domain(\comp(u)) = \domain(\fst(u))$,
(K4) $\codomain(\comp(u)) = \codomain(\snd(u))$ and
\begin{itemize}
\item[(K5)] $\fst(u) = \fst(v), \snd(u) = \snd(v) \Longrightarrow u=v$
\item[(K6)] $\domain(f) = \codomain(g) 
\Longrightarrow (\exists u:{\cal C}_2)\, (\snd(u)  = f\, \&\, \fst(u) = g)$
\end{itemize}

We introduce abbreviations: for mappings $f,g,h$ write $h \equiv g \circ f$ for 
$(\exists u \in {\cal C}_2)[\fst(u)=f\, \&\, \snd(u) = g \, \&\, \comp(u) = h],$
that is, the diagram 
$$\bfig
\Atriangle/<-`>`>/[\cdot`\cdot`\cdot;f`g`h]
\efig$$
is composable and commutes. Write $k \circ h \equiv g \circ f$  if
there is a mapping $m$ so that $m \equiv g \circ f$ and 
$m \equiv k \circ h$,
that is, the following diagram composes and commutes
$$\bfig
\square[\cdot`\cdot`\cdot`\cdot; f`h`g`k]
\efig
$$

In terms of these abbreviations we can express the
monoid laws: (K7) $f \equiv f \circ (\idmap_{\domain(f)})$,
(K8) $f \equiv (\idmap_{\codomain(f)}) \circ f$, and (K9)  if
 $k \equiv f \circ g$ and $\ell \equiv g \circ h$ then $k \circ h  \equiv f \circ \ell$.

\medskip
$f:a \to b$ and $a \to^f b$  are abbreviations for the conjunction
 $\domain\, f = a \,\&\, \codomain\, f = b$. We shall often omit
 $\circ$ and write $h\equiv gf$ for $h \equiv g \circ f$. Moreover
$\equiv$ is often replaced by $=$ when there is no danger of confusion.

\section{Subobjects and Relations}

We may define the notion of an $n$-ary relation in any category.
Recall that a mapping $f:A \to B$ is {\em monic} or {\em is a mono} if for any
mappings $h,k: U\to A$ with $fh = fk$ it holds that
$h=k$.  We write in this case $f: A \mono B$.  This notion can be
generalized to several mappings.  A sequence of mappings $r_1: R \to
X_1, \ldots, r_n: R \to X_n$ are {\em jointly monic,} if for any
$f,g:U \to R$
$$r_1 f = r_1 g, \ldots,
r_n f = r_n g \Longrightarrow f= g.$$
In this case we write $(r_1,\ldots,r_n): R \mono (X_1,\ldots,X_n)$.
We regard this as an {\em $n$-ary relation between the objects
$X_1,\ldots,X_n$.}  In particular, a {\em binary relation between
$X_1$ and $X_2$} is a pair of mappings $r_1: R \to X_1$ and $r_2: R \to
X_2$ which are jointly monic.  Another particular case is: if the
category has a terminal object $\terminal$, a $0$-ary relation $(): R
\mono ()$ means that the unique map $R \to \terminal$ is a mono.

Consider a category ${\cal C}$ with a terminal object $\terminal$. An
{\em element} of an object $A$ of ${\cal C}$ is a mapping
$x:\terminal \to A$. For a monic $m: M \to X$ and element $x$
of $X$ write $x \memof m$ if $(\exists a:\terminal \to M) ma
=x$. We say that {\em $x$ is a member of $m$.}
More generally, if $(m_1,\ldots, m_n): M \mono (X_1,\ldots,X_n)$
and $(x_1,\ldots, x_n) : \terminal \mono (X_1,\ldots,X_n)$ we write
$(x_1,\ldots,x_n) \memof (m_1,\ldots, m_n)$ if there is
$a:\terminal \to M$ so that $m_i a = x_i$ for all
$i=1,\ldots,n$.

To simplify notation we often write $x \in X$ and $(x_1,\ldots,x_n)
\in (X_1,\ldots,X_n)$ for $x : \terminal \to X$ and 
$(x_1,\ldots,x_n) : \terminal \to (X_1,\ldots,X_n)$, respectively. Note the
difference between the signs $\in$ (elementhood)  and $\memof$ (membership).

\medskip
We shall be interested in categories where there is no canonical
construction for products, but where it is merely assumed that they exist.
Recall that an {\em $n$-ary product diagram} in a category is a
sequence of mappings $X \to_{p_i} X_i$ ($i=1,\ldots,n$) so that for any
sequence of mappings $C \to_{p_i} X_i$ ($i=1,\ldots,n$) there is a unique
$h:C \to X$ such that $f_i \equiv h p_i$ for all $i=1,\ldots,n$.
We write
$$h \equiv \langle f_1,\ldots,f_n \rangle_{\bar{p}}$$
when $f_i \equiv h p_i$ for all $i=1,\ldots,n$, 
where $\bar{p}=p_1,\ldots,p_n$. 
It is convenient to drop the subscripts $\bar{p}$ when
the product diagrams are obvious from the context.

\begin{prop} \label{prop31}
Suppose that $X \to_{p_i} X_i$ ($i=1,\ldots,n$) is a
  product diagram. If $(r_1,\ldots,r_n) : R \to (X_1,\ldots,X_n)$, 
$r':R \to X$
  and $r' \equiv \langle r_1,\ldots,r_n \rangle_{\bar{p}}$, then
$r'$ is monic iff $(r_1,\ldots,r_n)$ are jointly monic.
Moreover, for $(x_1,\ldots,x_n) \in (X_1,\ldots,X_n)$, $x' \in X$
with $x' \equiv \langle x_1,\ldots,x_n \rangle_{\bar{p}}$, we have
$$x' \memof r' \Longleftrightarrow (x_1,\ldots,x_n) 
\memof (r_1,\ldots,r_n). \; \qed$$

\end{prop}


A binary relation $f=(\xi,\upsilon) : R \mono (X,Y)$ is a {\em partial
  function} in case $\xi$ is mono. It is a {\em total function} in
  case $\xi$ is iso. A relation $$f=(\xi_1,\ldots,\xi_n,\upsilon) : 
R \mono (X_1,\ldots,X_n,Y)$$ is a {\em partial function of $n$
  variables}
if $(\xi_1,\ldots,\xi_n) : R \mono (X_1,\ldots,X_n)$.  We write 
$$f : (X_1,\ldots,X_n) \pto Y.$$
It is {\em
  total
function of $n$ variables} if $R\to^{\xi} X_i$ ($i=1,\ldots,n$) is a
product diagram. We write $$f : (X_1,\ldots,X_n) \to Y.$$
For $x_1 \in X_1,\ldots, x_n \in X_n$ and $y \in Y$
we write 
$$f(x_1,\ldots,x_n) \equiv y$$
in case $(x_1,\ldots,x_n,y) \memof f$.

\section{Axioms of ETCS and CETCS}

Lawvere's theory ETCS (Lawvere 2005) has eight axioms: (L1) finite
roots exist, (L2) the exponential of any pair of objects exist, (L3)
there is a Dedekind-Peano object, (L4) the terminal object is
separating, (L5) axiom of choice, (L6) every object not isomorphic to
an initial object contains an element, (L7) Each element of a sum is a
member of one of its injections, (L8) there is an object with more than
one element.

We present a constructive version of ETCS, called CETCS, and some
extensions, by laying down axioms for a category $\cal C$.
(It should be evident that the following axioms may be formulated in
first-order logic in a language with ${\cal C}_0, {\cal C}_1, {\cal C}_2$
as sorts and the function symbols 
$\idmap, \domain, \codomain, \comp, \fst,\snd$ as indicated in
Section 2.)

Lawvere's (L1) says that the category is {\em bicartesian,} i.e.\ 
both cartesian and cocartesian.

Recall that $\cal C$ is {\em cartesian} if the conditions (C1) -- (C3)
are satisfied:

(C1) There is a terminal object $\terminal$ in  $\cal C$.

(C2) Binary products exist: For any pair of objects $A$ and $B$
there exists an object $P$ and two mappings   
$$A \to/<-/^p P \to^q B$$
which are such that if $A \to/<-/^f X \to^g B$ then there exists a
unique $h: X \to P$ so that $p h \equiv f$ and $q h \equiv g$.

(C3) Equalizers exist: For any parallel pair of mappings $A \two^f_g B$
there exists a mapping $e:E \to A$ so that $f e \equiv g e$
 and such that whenever $h:X \to A$ satisfies $f h \equiv
g h$ then there exists a 
unique $k: X \to E$ with $e k  \equiv h$.

A category $\cal C$ is {\em cocartesian} if it satisfies (D1) -- (D3),
which are the categorical duals of (C1) -- (C3).

(D1) There is an initial object $\initial$ in $\cal C$.

(D2) Binary sums exist: For any pair of objects $A,B$   there is
a diagram
\begin{equation} \label{sumdia}
A \to^i S \to/<-/^j B
\end{equation}
such that if $A \to^f T \to/<-/^g B$ then there is a unique $h :S \to T$ with
$h i \equiv f$ and $h j \equiv g$. 

(D3) Coequalizers exist: For any parallel pair of mappings $A \two^f_g B$
there exists a mapping $q:B \to Q$ so that $q  f  \equiv q g$
and such that whenever $h:B \to Y$ satisfies $h f \equiv h  g$ then
there exists a unique $k: Q \to E$ with $k q \equiv h$.

\medskip
The axiom (L2) of ETCS says together with (L1) that the category is
cartesian closed. Instead, we take for an axiom the following ($\Pi$) which,
together with cartesianess and axiom (G) below, states that the category
is locally cartesian closed. (This axiom is a theorem of ETCS.)

($\Pi$) Dependent products exist: For any 
mappings  $Y \to^g X \to^f I $  there exists a commutative diagram
\begin{equation}  \label{defpid}
\bfig                                                                          
\qtriangle/<-`->`/[Y``;\piev`g`]         
\square(550,0)[P`F`X`I;\pi_1`\pi_2`\vf`f]
\efig
\end{equation}
where the square is a pullback, and which is such that
for any element $i \in I$  and any partial function
$\psi=(\xi,\upsilon) : R \pto (X,Y)$
such that
\begin{itemize}
\item[(a)] for all $(x,y) \in (X,Y)$,  
$(x,y) \memof \psi$ implies $g y \equiv x$ and $f x \equiv i,$
\item[(b)] if $f x \equiv i$, then there is 
$y \in Y$ with $(x,y) \memof \psi$,
\end{itemize}
then there is a unique $s \in F$ so that $\vf s = i$ and 
for all $(x,y) \in (X,Y)$,  
\begin{equation} \label{equival1}
(s,x,y) \memof \alpha    \Longleftrightarrow 
(x,y) \memof \psi.
\end{equation}
Here $\alpha=(\pi_1,\pi_2,\piev): P \mono (F,X,Y)$.

A diagram (\ref{defpid}) satisfying these properties
 is called a {\em universal dependent product
diagram} or shortly a {\em universal $\Pi$-diagram} for $Y \to^g X \to^f I$.

\medskip
The third axiom (L3) of ETCS says, in now common terminology,
that there exists a {\em natural numbers
object} (NNO). A category $\cal C$ has an
NNO if there is a sequence of mappings (the NNO)
 $\terminal \to^0 N \to^s N$ so that for any other sequence of mappings
$\terminal \to^b A \to^h A$ there is a unique $f: N \to A$ with
$f 0 \equiv b$ and $f S \equiv h f$.

\medskip
Axiom (L4) states in modern terminology that $\terminal$ is a
separating object, i.e. as in Proposition \ref{prop41}. We consider instead
a stronger axiom (G) which is a theorem of ETCS.
A mapping $f:A \to B$ of $\cal C$ is {\em onto} if for any $y \in B$
there exists an $x \in A$ so that $y \equiv f x$. Our axiom is
 
\medskip
(G) Any mapping which is both onto and mono, is an isomorphism.

\medskip
The fifth axiom (L5) of ETCS states the axiom of choice in peculiar
way; see Section 5.2. A
more standard way is to first define an object 
$P$ of $\cal C$ to be a {\em choice object,} if for any onto $f:A
\to P$ there is a $g:P \to A$ with $f g = \idmap_P$.
The {\em axiom of choice} (AC) says that every object is a choice
object. This is a far too strong assumption in a constructive
setting. There is a constructively acceptable 
weakening which accords well with Bishop's distinction of operations
and functions, the {\em presentation axiom} (Aczel 1978):

\medskip
(PA) For any object $A$ there is an onto mapping $P \to A$ where $P$
is a choice object.

\medskip Axiom (L6) of ETCS says in contrapositive form: if an object has
no elements then it is an initial object. We take instead

\medskip
(I) The object $\initial$ has no elements.

\medskip
This together with (G) implies (L6).

\medskip
The Axiom (L7) of ETCS is {\em each element of a sum is a member of one of its 
injections.} We adopt this axiom unaltered but call it the {\em disjunction
principle (DP)} as it connects sums to disjunctions:

\medskip
(DP) In a sum diagram $A \to^i S \to/<-/^j B$:  
for any $z\in S$, $z \memof i$ or $z \memof j$.

\medskip
The final axiom (L8) of ETCS states that there exists object with at
least two elements. We state this as

\medskip
(NT, Non-triviality) For any sum diagram 
$\terminal \to^x S \to/<-/^y \terminal$ it
 holds that  $x \ne y$.

\medskip
There are two further axioms that we shall consider, which are
in fact theorems of ETCS.

\medskip
(Fct) Factorization. Any mapping $f$  can be factored as $f \equiv i e$
where $i$ is mono and $e$ is onto.

\medskip
(Eff) All equivalence relations are effective. For each equivalence relation $(r_1,r_2) : R \mono (X,X)$ there
is some mapping $e: X \to E$ so that 
$$(x_1,x_2) \memof (r_1,r_2) \Longleftrightarrow e x_1 \equiv e
x_2$$
for all $(x_1,x_2) \in (X,X)$.

\medskip
In summary, the theory CETCS consists of the axioms (C1 -- C3), (D1- D3), ($\Pi$),
(G), (PA), (I), (DP), (NT), (Fct) and (Eff). Observe that it is a finitely
axiomatized theory just as ETCS. We do not know whether this set of
axioms is optimal.

\begin{remarks} {\em
Note that it is not assumed that the (co)products or (co)equalizers
are given as functions of their data.  
The axiom (G) is in the terminology of Johnstone (2002)
that {\em $\terminal$ generates $\cal C$.} It entails that one can 
``reason using elements'' as the two following results exemplify. This
gives a substantial simplification of the internal logic.}
\end{remarks}

\begin{prop} \label{prop41}
Let $\cal C$ be a cartesian category which satisfies (G). Then 
\begin{itemize}
\item[(a)] For any pair of mappings $f,g: A \to B$, $f=g$ whenever 
$(\forall x\in A)(f x = g x)$.
\item[(b)] A mapping $f:A \to B$ is monic if and only if
$(\forall x,y\in A)(f x = f y \Rightarrow x = y)$.
\end{itemize}
\end{prop}
{\flushleft \bf Proof.} (b) follows easily from (a). To prove the non-trivial direction 
of (a):  assume that $(\forall x\in A)(f x = g x)$. Construct 
an equalizer $E \to^e A \two^f_g B$  of $f$ and $g$. Then $e$  
is monic. By the assumption and the equalizing property it is also
easy to see it is onto. Hence by (G) $e$ is an isomorphism. Since $f e = g e$ we get $f=g$. $\qed$

\medskip
Define an element-wise
inclusion relation for monos $m:M \to X $  and $n:N \to X$ 
$$m \subseteqof n \Longleftrightarrow_{\rm def} 
(\forall x \in X)(x \memof m \Rightarrow x \memof n)$$
The standard inclusion relation in a category is given by 
$m \le n\Longleftrightarrow_{\rm def} (\exists f:M \to N)(m = n f)$. Compare diagram (\ref{bishincl}).
Their correspondence is given by:

\begin{prop}  \label{elts2}
Let $\cal C$ be a cartesian category which satisfies (G). Then for all
 monos $m:M \to X $  and $n:N \to X,$
$$m \subseteqof n \Longleftrightarrow m \le n.$$
\end{prop}
{\flushleft \bf Proof.} ($\Leftarrow$) This is straightforward in any
category with a terminal object. ($\Rightarrow$) 
Suppose that $m:M \mono X $  and $n:N \mono X$ satisfies
$m \subseteqof n$. Form a pullback square 
$$\bfig \square[P`N`M`X;p`q`n`m ] \efig$$ 
To prove $m \le n$ it is evidently enough to show that $q$ is an
isomorphism.  Now $q$ is the pullback of a mono, so it is a mono as
well. By (G) it is sufficient to show that $q$ is onto. Let $y \in
M$. 
Thus $m  y \memof m$ and by assumption also $m
 y \memof n$. There is thus $t\in N$ with $m y
= n t$. Hence by the pullback square there is a unique
$u\in P$ so that $q u =y$ and $p u = t$. In
particular, this shows that $q$ is onto. $\qed$

\medskip
Functions as a graphs and as morphisms can be characterized as follows.

\begin{prop} \label{fungraph}
Let $\cal C$ be a cartesian category which satisfies (G). Let $r=(r_1,r_2):R \to (X,Y)$ be a relation.
Then
\begin{itemize}
\item[(a)] $r$ is a partial function if and only if 
\begin{equation}\label{pfcon}
(\forall x\in X)(\forall y,z \in Y)[(x,y) \memof r\; \&\; (x,z) \memof r \Rightarrow y=z].
\end{equation}
\item[(b)] $r$ is a total function if and only if 
\begin{equation}\label{pfcon2}
(\forall x \in X)(\exists ! y \in Y)(x,y) \memof r.
\end{equation}
\item[(c)] {\em (Unique Choice)} If $(\forall x \in X)(\exists ! y \in Y)(x,y) \memof r$, then
there is $f:X \to Y$ with 
$$(\forall x \in X)(x,f x) \memof r.$$
\end{itemize}
\end{prop}
{\flushleft \bf Proof.} (a): by definition $r$ is a partial function if and only if $r_1$ is mono.
By Proposition \ref{prop41}, $r$ is thus a partial function precisely when
$$(\forall s,t \in R)[r_1 s= r_1 t \Rightarrow s=t].$$
This is easily seen to be equivalent to (\ref{pfcon}). 

(b, $\Rightarrow$): Suppose $r$ is a total function. Then $r_1$ is iso. For $x \in X$, we have $(x,y) \memof r$ with $y= r_2 r^{-1}_1 x)$.
By (a) it follows that $y$ is unique.

(b, $\Leftarrow$):  Suppose (\ref{pfcon2}) holds. By (a) $r_1$ is mono.  For each  $x \in X$ there is some $t \in R$ and $y=r_2 t$ so that $(x,y) \memof r$. Thus $r_1$ is onto, and by (G) $r_1$ is iso.

(c): This is clear from (b, $\Leftarrow$) since then $r_1$ is invertible, and we may take $f=r_2 r^{-1}_1$:  for $x \in X$, $x= r_1 r^{-1}_1 x$ and $f x = r_2 r^{-1}_1 x$ so $(x,fx) \memof r$. $\qed$

\begin{prop} Let $\cal C$ be a cartesian category which satisfies (G). Then a commutative
diagram
$$\bfig
\square[P`A`B`C;\pi_1`\pi_2`f`g]
\efig$$
is a pullback diagram if, and only if, 
\begin{equation} \label{pbchar}
(\forall x\in A)(\forall y \in B)[fx = gy \Longrightarrow (\exists ! t \in P) x = \pi_1 t \; \&\; y = \pi_2 t].
\end{equation}
\end{prop}
{\flushleft \bf Proof.} ($\Rightarrow$) Immediate. ($\Leftarrow$): Assume (\ref{pbchar}). It follows that
$\pi_1$ and $\pi_2$ are jointly monic. Suppose there is given a commutative square
$$\bfig
\square[Q`A`B`C;q_1`q_2`f`g].
\efig$$
Form the pullback
$$\bfig
\square[R`Q`P`A\times B;h`k`\langle q_1,q_2 \rangle`\langle \pi_1,\pi_2 \rangle].
\efig$$
Clearly $h$ is mono, since it is a pullback of a mono. By (\ref{pbchar})
\begin{equation} 
(\forall s \in Q)(\exists ! t \in P) [q_1 s= \pi_1 t \; \&\; q_2 s = \pi_2 t].
\end{equation}
but this implies that $h$ is onto. Hence $h$ is iso by (G). Thus $m=kh^{-1}: Q \to P$ satisfies
$\pi_i m =q_i$ for $i=1,2$, and is the desired map.  It is unique since $\pi_1$ and $\pi_2$ are joint monic. $\qed$

\section{Basic Set-theoretic Consequences}
\medskip
We mention some easy consequences of the axioms.

\begin{prop}  \label{qprop}
{\em (Quotient sets.)} Suppose that the bicartesian 
category $\cal C$ satisfies (G).
For any equivalence relation $r=_{\rm def } (r_1,r_2) : R \mono (X,X)$ there
exists a mapping
$q: X \to Q$ so that for all $(x_1,x_2) \in (X,X)$
\begin{equation} \label{rqequiv}
(x_1,x_2) \memof r \Longrightarrow q x_1 = q x_2
\end{equation}
and if $f: X \to Y$ is any mapping with
\begin{equation} \label{rqequiv2}
(x_1,x_2) \memof r \Longrightarrow f x_1 = f x_2.
\end{equation}
then there exists a unique $h:Q\to Y$ with $h q = f$.

In case the category also satisfies (Eff) it follows that
(\ref{rqequiv}) is an equivalence.
\end{prop}
{\flushleft \bf Proof.} Construct a coequalizer diagram
$$R \two^{r_1}_{r_2} X \to^q Q.$$
Since the diagram commutes, the implication (\ref{rqequiv}) holds.
Let $f:X \to Y$ be any mapping satisfying the implication
(\ref{rqequiv2}). Thus for any $t \in R$, $f r_1  t = 
f  r_2  t$. Thus by Proposition \ref{prop41} (a) we have
$f r_1 = f  r_2$ and since $q$ is a coequalizer, there
is a unique $h:Q\to Y$ with $h q = f$.

From Axiom (Eff) it follows that there is some $e: X \to E$  such that
\begin{equation} \label{rqeuiv3}
(x_1,x_2) \memof r \Longleftrightarrow e x_1 = e
x_2
\end{equation}
for all $(x_1,x_2) \in (X,X)$. Thus $e r_1 = e r_2$.
Let $e': Q \to E$ be the unique mapping so that $e' q = e$.
Thus if $q x_1 = q x_2$, it follows that
 $e x_1 = e x_2$ and
hence $(x_1, x_2) \memof r$ by (\ref{rqeuiv3}). $\qed$

\begin{prop} {\rm (Induction.)} Assume that $\cal C$ is a cartesian category
  which satisfies (G) and (NNO). Let $r: R \mono N$. 
Suppose that $0 \memof r$ and that
for each $n\in N$, $n \memof r$ implies $S n \memof r$. Then for
all $n \in N$, $n \memof r$.
\end{prop}
{\flushleft \bf Proof.} Since $0 \memof r$, there is $z: \terminal
 \to R$ with $0 \equiv r z$. Form a pullback square
$$\bfig
\square[P`R`R`N;q`p`r`S \circ r]
\efig
$$
As $r$ is mono, so is $p$. We claim that $p$ is onto. Let $u : \terminal
\to R$. Thus $r u \memof r$. Hence by assumption $ S r
u \memof r$. There is thus a map $v: \terminal \to R$ with $S r
u  = r v$. By the pullback property there is $x: \terminal
\to P$ so that $p x = u$ and $q x = v$. In particular $p$
is onto. By (G) $p$  is an isomorphism. Let $p^{-1}$ be its inverse.
Thus $q p^{-1} : R \to R$. By the property of the natural
number object there is a unique $f:N \to R$ with $f 0 = z$
and $f S = q p^{-1} f$. Now $(r  f) 0 = 0$
and $$(r f) S = r  q  p^{-1} f 
= S r f.$$
But $\idmap_N$ instead of $r \circ f$ also satisfies these two
equations. Thus $r f = \idmap$. Thus for any $n\in N$, $r
f n = n$, and hence $n \memof r$. $\qed$

\begin{prop} {\em (Exponential objects.)} Assume that $\cal C$ is a
  cartesian category that satisfies (G) and ($\Pi$).
Then for any objects $X$ and $Y$ there is an object $E$ and a
total function $e: (E,X) \to Y$ such that for every morphism
$f: X \to Y$ there is a unique $s \in E$ such that 
for $x \in X$ and $y \in Y$:
$$e(s,x) \equiv y \Longleftrightarrow f \circ x \equiv y$$
\end{prop}

\begin{thm} {\em (Dependent choices.)} Assume that { }$\cal C$ is a cartesian
category that satisfies (G), ($\Pi$), (Fct) and (PA). 
Then for any object $X$, any total relation
     $r=(r_1,r_2):R \mono (X,X)$ and any $x \in X$ there is a
morphism $f:N \to X$ with $f 0 = x$ and for all $n \in N$
\begin{equation} \label{dcrelation}
(f n, f \circ S n) \memof r.
\end{equation}
\end{thm}
{\flushleft \bf Proof.} (Sketch) Take a projective cover  $p:P \to X$
of $X$. Since $r$ is total, we have thus for each $u \in P$ some $v
\in P$ with $(p u, p v) \memof r$.
As $P$ is a choice object, there is a morphism $g: P \to P$ with 
$(p u, p g u) \memof r$ for all $u \in P$. 
Let $x \in X$. Then $p \circ w \equiv x$ for some $w \in
P$. Now $\terminal \to^0 N \to^S N$ is a natural numbers object, so
there is $h:N \to P$ with $h 0 = w$ and $h S = g h$. Now it is easy to
check by induction that $f =_{\rm def} p h$ satisfies 
(\ref{dcrelation}). $\qed$

\subsection{Constructing New Relations}

We review some of the possibilities to construct relations in a
 bicartesian category satisfying the axioms (G), ($\Pi$), (DP), (Fct)
 and (I).
On any object $X$ the identity mapping gives a universally true relation 
${\sf t}_X = \idmap_X : X \to X$, i.e.\ for all $x\in X$
$$x \memof {\sf t}_X.$$  
The unique mapping from the initial object ${\sf f}_X: \initial \to X$ gives
an universally false relation, i.e.\ for all $x \in X$,
$$\lnot (x \memof {\sf f}_X).$$
If $E \to^e X \two^g_h Y$ is an equalizer diagram, then for $x \in X$
$$x \memof e \Longleftrightarrow g x = h x.$$
Given a relation $r=(r_1,\ldots,r_n): R \mono (X_1,\ldots,X_n)$ we can
extend it with a variable. Let $Y$ be a object and let $R \to/<-/^p
R'\to^q Y$ be a product diagram. The extended relation
$$r'= (r_1 p,\ldots,r_n p,q) : R' \mono (X_1,\ldots,X_n, Y)$$
satisfies, for all $(x_1,\ldots,x_n,y) \in (X_1,\ldots,X_n,Y)$ that
$$(x_1,\ldots,x_n,y) \memof r' \Longleftrightarrow (x_1,\ldots,x_n) \memof r.$$
If $\sigma: \{1,\ldots,n\}\to \{1,\ldots,n\}$ is a permutation then
$$r_{\sigma}= (r_{\sigma(1)},\ldots,r_{\sigma(n)}): R \mono 
(X_{\sigma(1)},\ldots,X_{\sigma(n)})$$ satisfies
for all $(x_1,\ldots,x_n) \in (X_{\sigma(1)},\ldots,X_{\sigma(n)})$
$$(x_1,\ldots,x_n) \in r_{\sigma}
\Longleftrightarrow 
(x_{\sigma^{-1}(1)},\ldots,x_{\sigma^{-1}(n)}) \in r.
$$

The following lemma is standard
\begin{lemma} If in the universal $\Pi$-diagram
\begin{equation} 
\bfig                                                                          
\qtriangle/<-`->`/[Y``;\piev`g`]         
\square(550,0)[P`F`X`I;\pi_1`\pi_2`\vf`f]
\efig
\end{equation}
the mapping $g$ is mono, then so is $\vf$. $\qed$
\end{lemma}

Relations can be combined using the logical operations ($\land$,
$\lor$, $\Rightarrow$) and
quantifiers ($\forall$,$\exists$) over fixed objects:

\begin{thm} \label{elemcomp}
Let $\cal C$ be a bicartesian category satisfying the axioms
 (G), ($\Pi$), (DP), (Fct) and (I).
Let $r=(r_1,\ldots,r_n): R \mono (X_1,\ldots,X_n)$ and 
$s=(s_1,\ldots,s_n): S \mono (X_1,\ldots,X_n)$.
Then exists $(r \land s),(r \lor s), (r \Rightarrow s): 
R \mono (X_1,\ldots,X_n)$ so that for all
$x=(x_1,\ldots,x_n)\in (X_1,\ldots,X_n)$
\begin{itemize}
\item[(a)] $x\memof (r \land s)$ if and only if $x \memof r$ and $x \memof s$,
\item[(b)] $x\memof (r \lor s)$ if and only if $x \memof r$ or $x \memof s$,
\item[(c)] $x\memof (r \Rightarrow s)$ if and only if $x \memof r$ implies $x \memof s$,
\end{itemize}
Moreover, if  $m: M \mono (X_1,\ldots,X_n,Y)$ then there is 
$\forall(m): A \mono (X_1,\ldots,X_n)$ and $\exists(m): E \mono (X_1,\ldots,X_n)$
so that for all
$x=(x_1,\ldots,x_n) \in (X_1,\ldots,X_n)$
\begin{itemize}
\item[(d)] $x\memof \forall(m)$ if and only if for all $y\in Y$,
$(x_1,\ldots,x_n,y) \memof m$,
\item[(e)] $x\memof \exists(m)$ if and only if for some $y\in Y$,
$(x_1,\ldots,x_n,y) \memof m$.
\end{itemize}
\end{thm}
{\flushleft \bf Proof.} By Proposition \ref{prop31} it is enough to
prove (a) -- (c) for the case when $n=1$; write $X=X_1$, $r=r_1$,
$s=s_1$. 

As for (a): form the
pullback square
$$\bfig\square[P`S`R`X;q`p`s`r]\efig.$$
The diagonal, call it  $(r \land s)$ is a mono. It is straightforward
by the pullback property
 that the equivalence in (a) holds.

As for (b): form a sum diagram $R \to^i U \to/<-/^j S$. Let $f: U \to
X$ be the unique mapping with $r = f i$ and $s = f j$. Let
$U \to^e I \to^m X$ be a factorization of $f$ as an onto mapping
followed by a mono (Fct).  We claim that  $(r\lor
s)=_{\rm def} m$ satisfies the equivalence in (b).
Suppose that $x \in X$ satisfies $x \memof r$. Then $x = r t$ for
some $t \in R$. Thus $x = f i t = m e i
t$, and hence $x \memof m$. Similarly $x \memof s$ implies $x \memof
m$. Suppose on the other hand that $x \memof m$. Now, $e$ is onto so 
there is some $u \in U$ with $x = f u$. Then by Axiom (DP) we
have $u= i t$ for some $t \in R$, in which case $x \memof r$, or we
have $u= j v$ for some $v \in S$, in which case $x \memof s$.

As for (c): Form the pullback
\begin{equation} \label{casec1}
\bfig
\square/ >->` >->` >->` >->/[Q`S`R`X;q`p`s`r]
\efig
\end{equation}
Axiom ($\Pi$) yields for $Q \to^p R \to^r X$ a universal $\Pi$-diagram
\begin{equation} \label{casec2}                             
\bfig
\qtriangle/<-`->`/[Q``;\piev`p`]
\square(550,0)[P`F`R`X;\pi_1`\pi_2`\vf`r]                                    
\efig                                         
\end{equation}
We claim that $(r \Rightarrow s) =_{\rm def} \vf$ makes the equivalence
in (c) true. Let $x \in X$. To prove ($\Rightarrow$)
assume that $x \memof \vf$. Thus $x = \vf u$ for some $u \in F$.
Suppose $x \memof r$. Thus $x= r v$ for some $v \in R$. By the 
pullback in (\ref{casec2}) there is $w:\terminal \to w$ so that $\pi_2
 w = v$ and $\pi_1 w = u$. We have further by the diagrams
$$x= r v = r \pi_2 w = 
r p \,\piev\, w = 
s q \,\piev\, w.$$
This shows $x \memof s$. As for the converse ($\Leftarrow$) suppose
the implication
$$x \memof r \Rightarrow x \memof s$$
holds. We aim to show $x \memof \vf$ using the properties of the
universal $\Pi$-diagram. Form a pullback diagram
\begin{equation} \label{casec3}
\bfig
\square/ >->` >->` >->` >->/[T`Q`\terminal`X;t``r p`x]
\efig
\end{equation}
Then $\psi= (p t,t):T \mono (R,Q)$ is evidently a partial function
 since $p$ and $t$ are monic. If $\psi(u) \equiv v$, then
there is $w \in T$, so that $u= p  t w$ and $v=t w$,
and hence $u = p v$ and $r u = r p v =x$. This
verifies condition (a) of ($\Pi$). To verify condition (b) 
of ($\Pi$), assume that
$r u = x$. Thus $x \memof r$, and so by the implication above
$x \memof s$, i.e.\ $s w = x$, for some $w$. 
By the pullback (\ref{casec1})
there is $v: \terminal \to Q$ with $u = p v$ and $w = q
v$. Thus $r p v = x$. But, by the pullback (\ref{casec3}) there
is $z: \terminal \to T$ with $t z = v$. Now $(p t 
z,t z) = (u,v)$, i.e.\ $\psi(u) \equiv v$. According to ($\Pi$),
there is now some $k \in F$ with $\vf k =x$. Thus $x \memof \vf$.

As for (d): Suppose $m=(m_1,m_2): M \mono (X,Y)$. Using (Fct) factor
$m_1$ into a onto mapping followed by a mono $M \to^e I \to^i X$.
Let $\exists(m) = i$.
Thus using the fact that $e$ is onto
$$x \memof \exists(m) \Leftrightarrow
(\exists t\in I) x= i  t \Leftrightarrow
(\exists s \in M) x = i e  s \Leftrightarrow
(\exists s \in M) x = m_1  s.
$$
The latter implies that $(x, m_2 s) \memof
(m_1,m_2)$. Clearly $m_2 s \in Y$. Conversely,
suppose that for some $y \in Y$ we have $(x,y) \memof (m_1,m_2)$. Thus
for some $s \in M$ it holds that $x= m_1 s$ and $y= m_2 s$ and we 
have $x \memof \exists(m)$.

As for (e): Suppose $m=(m_1, m_2): M \mono (X,Y)$. First construct a
product diagram $X \to/<-/^p U \to^q Y$. Then let $m' \equiv \langle
m_1,m_2 \rangle_{p,q}$. Use ($\Pi$) to obtain a universal $\Pi$-diagram
\begin{equation} \label{casee}                             
\bfig
\qtriangle/<-`->`/[M``;\piev`m'`]
\square(550,0)[P`F`U`X;\pi_1`\pi_2`\vf`p]                                    
\efig                                         
\end{equation}
We let $\forall(m) = \vf$. Suppose $x \in X$. To prove (e, $\Rightarrow$)
suppose $x \memof \vf$ and $y \in Y$. Thus $x = \vf  f$ for some 
$f \in F$ and moreover there is $u \in U$ with $x= p  u$ and 
$y = q  u$. By the pullback in (\ref{casee}) we get $w \in P$ so
that $u = \pi_2  w$ and $f = \pi_1  w$. From the triangle of
(\ref{casee}) it follows that $ m' \,\piev\, w = 
\pi_2  w $. Hence $u \memof m'$ and thus $(x,y) \memof m$. 
To prove (e, $\Leftarrow$)
let $x \in X$ be fixed and suppose that for all $y \in Y$, $(x,y) 
\memof m$. Let $n: N \to M$ be the pullback of $x$ along $m_1$:
\begin{equation} \label{casee2}
\bfig\square[N`M`\terminal`X;n``m_1`x]\efig
\end{equation}
Then $(m' n, n) : N \mono (U,M)$ is a partial function since
both $m'$ and $n$ are mono.  As for  condition (a): if 
$(u,v) \memof  (m'  n, n)$ then $u = m'  n  t$ and 
$v= n  t$ for some $t \in N$. Clearly, $m'  v= u$ and $p
 u =m_1  n  t = x$. Regarding condition (b): Suppose
that $u\in U$ satisfies $p  u = x$. Let $y = q  u$. By the
first assumption $(x,y) \memof m$. Thus for some $s \in M$, $x= m_1
 s$ and $y = m_2  s$. By construction of $m'$ we have $m'
 s = u$. Since $x= m_1  s$, the pullback (\ref{casee2})
gives a unique $t\in N$ with $s = n  t$. Thus $t$ is a witness
to $(u,s) \memof (m'  n, n)$.
Since conditions (a) -- (b) are now verified,  ($\Pi$) gives $f \in F$
satisfying, in particular, $\vf  f = x$. Hence $x \memof \vf$. 
$\qed$

\subsection{Decidable Relations and Classical Logic}

Let $\cal C$ be a CETCS category. Construct a two element set using
the sum axiom  $\terminal \to^f \bool \to/<-/^t\terminal$.
If $r: P \mono X$ is decidable,
i.e.\ for all $x \in X$,
$$\mbox{$x \memof r$ or $\lnot x\memof r$},$$
then we can construct $\chi_r: X \to \bool$ so that for all $x \in X$
$$\mbox{$x \memof r \land \chi_r(x) = t$ or 
$(\lnot x\memof r) \land \chi_r(x) = f$},$$
It follows that $\chi_r$ is the unique map $X \to \bool$ such that $x \memof r$ iff
$\chi_r(x) =t$. Thus $\terminal \to^t \bool$ classifies decidable
relations. In case we take the axioms of CETCS with classical 
logic every relation is decidable, and hence $\terminal \to^t \bool$
is a full subobject
classifier for the category. In this case $\cal C$ is a topos.

\medskip
The Lawvere's choice axiom (L5) states: If $f:A \to B$ is mapping and
$A$ contains at least one element, then there is a mapping $g:B \to A$ so that 
$fgf=f$.

\begin{thm} In CETCS with classical logic (AC) and (L5) are equivalent.
\end{thm}

\begin{cor} ETCS and CETCS +PEM + AC have the same theorems. 
\end{cor}

\section{Correspondence to Standard Categorical Formulations}

\medskip
\begin{lemma} \label{lm61}
Let $\cal C$ be a cartesian category which satisfies (G). Then
a pullback of an onto mapping is again an onto mapping. 
\end{lemma}

We recall some basic definitions from (Johnstone 2002):

\begin{defin}
A sequence of mappings $A \to^e I \to^m B$ is an {\em image factorization} of 
$f:A \to B$ if  $f \equiv m \circ e$ and $m$ is a mono, 
and whenever $f \equiv m' \circ e'$
where $m':I' \to B$ is mono then there is some $t:I \to I'$ with $m
\equiv m' \circ 
t$. Such an $m$ is a called an {\em image of $f$}.
\end{defin}

\begin{defin} A morphism $f$ is a {\em cover} if 
whenever it can be factored as $f \equiv m \circ g$ where
$m$ mono, then $m$ be must an isomorphism.
\end{defin}

\begin{prop} \label{imfactcov}
In any category, if $A \to^e I \mono^i B$ is an image factorization of $f:A \to B$, then $e$ is a cover.
\end{prop}

\begin{thm} \label{ontomono}
Let $\cal C$ be a cartesian category satisfying (G). If $f: X \to Y$ is factored as $X \to^g I \mono^i Y$, where 
$g$ is onto, then it is an image factorization.
\end{thm}
{\flushleft \bf Proof.} Suppose therefore that $X \to^h J \mono^j Y$
is another factoring of $f$. It is sufficient to show that $i \le j$
as subobjects of $Y$. By Proposition \ref{elts2} it is equivalent to
prove $i \subseteqof j$. Suppose that $y \in Y$ satisfies $y \memof i
$. Then $y= i t$ for some $t \in I$. Now $g$ is onto, so there is $x
\in X$ with $g x = t$. Now $y= i g x = f x = j (h x)$.  Hence $y
\memof j$. Thus we have $i \subseteqof j$. $\qed$

\begin{lemma} \label{lm66}
Suppose that $\cal C$ is a cartesian category that 
satisfies (G). Then 
\begin{itemize}
\item[(a)] every onto mapping is a cover,
\item[(b)] if $\cal C$ in addition satisfies (Fct), then
 every cover is onto.
\end{itemize}
\end{lemma}
{\flushleft \bf Proof.} (a): If $f:A\to B$ is onto, then $A\to^f B \to
^{\idmap}B$ is an image factorization, so by Theorem \ref{ontomono} and 
Proposition \ref{imfactcov} 
$f$ is a cover.

(b): Let $f:A\to B$ be a cover. By (Fct) take a factorization $A\to^e I
\to^i B$ of $f$ where $e$ is onto and $i$ is mono. Now since $f$ is a
cover, $i$ is an isomorphism. Hence $f$ is onto as well. $\qed$

\medskip
In standard category-theoretic terms (Johnstone 2002) various combinations of
the CETCS axioms can be characterized by the following theorems. First
recall that a {\em regular category} is a category with finite limits,
which has image factorization and where covers are preserved by
pullbacks. 

\begin{thm} \label{regchar}
Let $\cal C$ be a cartesian category satisfying (G). Then $\cal C$ 
satisfies (Fct) if,
and only if, $\cal C$ is a regular category 
where the terminal object is projective.
\end{thm}
{\flushleft \bf Proof.} ($\Rightarrow$) According 
to Theorem \ref{ontomono} $\cal C$ has
image factorizations. By Lemma \ref{lm66}   onto morphisms are 
the same as covers. Thus by Lemma \ref{lm61}
covers are preserved by pullbacks. This shows that $\cal C$ is
regular.
If $A \to \terminal$ is a cover then it is onto, and hence
$\terminal$ is a choice object. Since the category is regular, it
follows that $\terminal$ is projective.

($\Leftarrow$) 
Suppose $\terminal$ is projective. Hence any cover is
onto. 
Thus by regularity, any morphism can be factored as an onto morphism
followed by a mono. This gives (Fct). $\qed$

\begin{thm} \label{lccchar} 
Let $\cal C$ be a cartesian category satisfying (G). Then
$\cal C$ is locally cartesian closed if and only if $\cal C$ satisfies
the axiom ($\Pi$).
\end{thm}
{\flushleft \bf Proof.} See Section 7. $\qed$

\begin{lemma} \label{balance}
In a CETCS category $\cal C$ every epi is onto; consequently $\cal C$ is balanced. 
\end{lemma}
{\flushleft \bf Proof.} Let $f:A \to B$ be an epi. Form the sum
$\terminal \to^i S \to/<-/^j B$. Let $m: M \mono B$ be a subobject
so that $y \memof m$ iff $(\exists x\in A) f x = y$. Then form the sum
$\terminal \to^r K \to/<-/^s M$ and let $k: K \to S$ be the unique
mapping so that $k r = i$ and $k s = j m$. Define, using Section 5.1,
 an equivalence
relation $(r_1,r_2) : R \to (S,S)$ by 
$$(u,z) \memof (r_1,r_2)  \Longleftrightarrow ((\exists w \in K) kw=u \Leftrightarrow
(\exists w \in K) kw=z).
$$
By Proposition \ref{qprop} let $q: S \to Q$ be such that
$$(u,z) \memof (r_1,r_2)  \Longleftrightarrow  qu = qz.$$
Let $g:B \to S$ be given by $g = i \circ !_B$ and $h=j:B \to S$.
It is straightforward to check that for all $x\in A$, 
$qgfx = qhfx$. Thus $qgf=qhf$, and since $f$ is epi, $qg=qh$.
For each $y\in B$ we have,  since $(\exists w \in K) kw=gy$ is true,
that 
$$(\exists w \in K) kw = hy.$$
By (DP) and disjointness of sums we must have $w = st$ for some $t \in M$.
Hence $jmt = kst = kw=hy = jy$. Since $j$ is mono, $mt=y$. Thus $y \memof m$,
that is $(\exists x \in A) fx=y$.

The last statement follows by Axiom (G).
 $\qed$

\begin{thm} Let $\cal C$ be a category.
 Then $\cal C$ satisfies CETCS if,
and only if, $\cal C$ has the following properties 
\begin{itemize}
\item[(i)] it is locally cartesian closed, 
\item[(ii)] it is a pretopos,
\item[(iii)] it has NNO,
\item[(iv)] its terminal object is projective and generates $\cal C$,
\item[(v)] $\initial \not\cong \terminal$,
\item[(vi)] it satisfies the disjunction property,
\item[(vii)] it has enough projectives.
\end{itemize}
\end{thm}
{\flushleft \bf Proof.} ($\Rightarrow$): (i) follows from Theorem
\ref{lccchar}. Properties (iii),(vi),(vii) are axioms of CETCS. (iv) follows
from Theorem \ref{regchar}. (v) is clear by (I). By (Lemma 1.5.13
--14, Johnstone 2002) every locally cartesian closed which is
cocartesian and balanced (Lemma \ref{balance}) is a pretopos.

($\Leftarrow$) It is known that in a locally cartesian closed pretopos
with NNO has all coequalizers (Remark 2.8, Moerdijk and Palmgren
2000). Using Theorem \ref{lccchar} we get axiom $(\Pi)$. Axioms (G),
(NNO), (PA) and (DP) are given. (I) follows easily from (v) using
uniqueness of mappings. In a pretopos the pullback object of $x$ and
$y$ in a sum diagram $\terminal \to^x S \to/<-/^y \terminal$ 
will be $\initial$, so (NT) follows from (I). In pretopos every
map can be factored as a cover followed by a mono. But using that
$\terminal$ is projective we can show that covers are onto, so (Fct)
is verified.  In a pretopos all equivalence relations are effective,
so (Eff) follows. $\qed$

\section{Functor-free Formulation of LCCCs} 

The standard way (Johnstone 2002) of defining a locally cartesian
category $\cal C$ is to say that it is a cartesian category so that
pullbacks along a mapping $f:X \to Y$ induces a functor $f^*: {\cal
C}/Y \to {\cal C}/X$ and that this functor has a right adjoint $\Pi_f
: {\cal C}/X \to {\cal C}/Y$. These functors must, in particular, be
defined on the objects of the slice categories. This means that the
pullback object must be possible to construct as a function of
mappings $g: A \to Y$ and $f:X \to Y$. This can be forced if one assumes the
full axiom of choice in the meta-theory of $\cal C$, but is not
possible if we only use intuitionistic logic. Makkai (1996) has
developed a theory of functors --- anafunctors --- by which one can avoid
such uses of choice. In (Palmgren 2008) we showed how LCCCs could be
formulated replacing $f^*$ and $\Pi_f$ by the appropriate anafunctors,
so that $\Pi_f$ is the right adjoint of $f^*$. We here extract what is
the existence condition for such $\Pi_f$ and formulate it without
functors. Thus a functor-free formulation of LCCC will be given in
Definition \ref{funfree}.

{\em A $\Pi$-diagram} for $Y \to^g X \to^f I$ is a commutative
diagram of the form

\begin{equation} \label{pidiag}
\bfig 
\qtriangle/<-`->`/[Y``;\piev`g`]
\square(550,0)[P`F`X`I;\pi_1`\pi_2`\vf`f] 
\efig                                 
\end{equation}  
where the square on the right is a pullback diagram.
The object $F$ is called the {\em parameter object} of the diagram.

If we have a second $\Pi$-diagram for $Y \to^g X \to^f I$

\begin{equation} \label{pidiag1b}
\bfig
\qtriangle/<-`->`/[Y``;\piev'`g`]
\square(550,0)[P'`F'`X`I;\pi'_1`\pi'_2`\vf'`f]
\efig                                 
\end{equation}  
we say that a mapping $t: F' \to F$ is a {\em $\Pi$-diagram morphism} from the second diagram
to the first diagram  if $\vf   t \equiv \vf'$ and the unique map $s:P' \to P$
such that $\pi_2  s \equiv \pi'_2$ and $\pi_1  s \equiv 
t  \pi'_1$ also satisfies $\piev  s \equiv \piev'$.

\begin{equation} \label{pidm}
\xymatrix{
&& {P'} \ar@{->}[rr]\ar@{->}[dll] \ar@{->}[ddr] \ar@{.>}[drr] && 
{F'}\ar@{->}[ddr] \ar@{.>}[drr] && \\
{Y}\ar@{->}[drrr]&&&& {P} \ar@{->}[rr] \ar@{->}[llll] \ar@{->}[dl] && {F}\ar@{->}[dl] \\
&&& {X} \ar@{->}[rr] && {I} & \\
}
\end{equation}

It is easily seen that the $\Pi$-diagrams and $\Pi$-diagram morphisms 
over fixed mappings
$Y \to^g X \to^f I$ forms a category.

A {\em universal $\Pi$-diagram for $Y \to^g X \to^f I$:} is a  
$\Pi$-diagram
\begin{equation} \label{pidiag2}
\bfig
\qtriangle/<-`->`/[Y``;\piev`g`]
\square(550,0)[P`F`X`I;\pi_1`\pi_2`\vf`f]
\efig 
\end{equation}
which is such that for any other $\Pi$-diagram  
\begin{equation} \label{pidiag3} 
\bfig 
\qtriangle/<-`->`/[Y``;\piev'`g`]
\square(550,0)[P'`F'`X`I;\pi'_1`\pi'_2`\vf'`f]
\efig
\end{equation}  
there is a unique mapping $n: F' \to F$ so that $\vf' \equiv \vf  n$ and
that the unique mapping $m: P' \to P$,  with $n \pi'_1 \equiv \pi_1  m$
and $\pi'_2 \equiv \pi_2  m$, satisfies $\piev' \equiv \,\piev\,  m$.

\begin{defin} \label{funfree}
{\em 
A cartesian category is {\em locally cartesian closed,} if it satisfies the
{\em generalized exponential axiom} or the {\em $\Pi$-axiom}: for every composable
pair of maps $Y \to^g X \to^f I$ there is an universal exponential diagram
as in (\ref{pidiag2}). That is, the category of $\Pi$-diagrams over 
$Y \to^g X \to^f I$ has a terminal object.}
\end{defin}

\subsection{Characterization of Universal $\Pi$-diagrams}

We have the following characterization of $\Pi$-diagrams where 
the parameter object is $F=1$.

\begin{lemma} \label{lem21}
Consider a cartesian category satisfying (G). Let $Y \to^g X \to^f I$ be morphisms and let $i \in I$ be an element.
For a pair of morphisms $\psi=(r_1,r_2):R \to (X,Y)$ the diagram
\begin{equation}  \label{defpid01}
\bfig                                                                          
\qtriangle/<-`->`/[Y``;r_2`g`]         
\square(550,0)[R`1`X`I;`r_1`i`f]
\efig
\end{equation}
is a $\Pi$-diagram if and only if 
\begin{itemize}
\item[(A1)] $\psi$ is a partial function (i.e. $r_1$ is mono)
\item[(A2)] $(\forall x \in X)[f x = i \Longrightarrow (\exists y \in Y)(x,y) \memof \psi]$
\item[(A3)] $(\forall x \in X)(\forall y \in Y)[(x,y) \memof \psi \Longrightarrow f x =i \land g y =x]$
\end{itemize}
\end{lemma}
{\flushleft \bf Proof.} ($\Rightarrow$) Suppose (\ref{defpid01}) is a
$\Pi$-diagram. Since $i$ is mono, the pullback diagram entails that
$r_1$ is mono. Hence $\psi$ is a partial function. Property (A2)
follows by the pullback property. (A3) follows since the whole diagram
is commutative.

($\Leftarrow$) Suppose that (A1) -- (A3) are satisfied. By (A3) it follows that the entire diagram commutes. (A1) and (A2) together yields that the square is a pullback. $\qed$

\begin{lemma} \label{lem22}
Consider two $\Pi$-diagrams in a cartesian category satisfying (G).
\medskip
\begin{equation}  \label{defpid02}
\bfig                                                                          
\qtriangle/<-`->`/[Y``;\piev'`g`]         
\square(550,0)[P'`F'`X`I;\pi_1'`\pi_2'`\vf'`f]
\qtriangle(1500,0)/<-`->`/[Y``;\piev`g`]         
\square(2050,0)[P`F`X`I;\pi_1`\pi_2`\vf`f]
\efig
\end{equation}
Let $\chi : F' \to F$ be such that $\vf \chi =\vf'$. There is a unique $\kappa: P' \to P$ so that $\pi_1 \kappa = \chi \pi'_1$ and 
$\pi_2 \kappa = \pi'_2$. For this $\kappa$ it holds that $\piev \kappa = \piev'$ if and only if for all $v \in F'$, $x \in X$ and $y \in Y$ 
$$(v,x,y) \memof (\pi'_1,\pi'_2,\piev')  \Longleftrightarrow (\chi v, x, y) \memof (\pi_1, \pi_2, \piev).$$
\end{lemma}
{\bf \flushleft Proof.} ($\Leftarrow$): Assume the equivalence. Let $t \in P'$ be arbitrary. We prove 
$\piev \kappa t = \piev' t$. Clearly $(\pi'_1 t,\pi'_2 t,\piev' t) \memof (\pi'_1,\pi'_2,\piev') $, so by  the equivalence $(\chi \pi'_1 t,\pi'_2 t,\piev' t)  \memof  (\pi_1, \pi_2, \piev)$. Thus  there is a $u \in P$ with $\chi \pi'_1 t = \pi_1 u$, $\pi'_2 t = \pi_2 u$ and  $\piev' t= \piev u$. 
Now we have $\pi_1 u = \chi \pi'_1 t =  \pi_1 \kappa t $ and $\pi_2 u = \pi'_2 t = \pi_2 \kappa t$. By the pullback property, $\pi_1$ and $\pi_2$ are jointly mono, so $u= \kappa t$. Thus $\piev' t= \piev \kappa t$.

($\Rightarrow$): Assume $\piev \kappa=  \piev'$.   Suppose $(v,x,y) \memof (\pi'_1,\pi'_2,\piev') $. Thus for some $t \in P'$, it holds that 
$v= \pi'_1 t$, $x=\pi'_2 t$ and $y=\piev'  t$. Hence $x= \pi_2 \kappa t$, $y= \piev \kappa t$. Finally $\pi_1 \kappa =\chi \pi'_1$ gives
$\chi v = \pi_1 \kappa t$, so that $(\chi v, x, y) \memof (\pi_1, \pi_2, \piev)$. For the converse, assume  $(\chi v, x, y) \memof (\pi_1, \pi_2, \piev)$.
Thus $\chi v= \pi_1 s$, $x=\pi_2 s$ and $y=\piev s$ for some $s \in P$. Then
$$f \pi_2 s = \vf \pi_1 s = \vf \chi v = \vf' v.$$
Thus there is a unique $t \in P'$ with $\pi'_2 t = \pi_2 s$ and $\pi'_1
t=v$. We have then $\pi'_2 t=x$, so to prove $(v,x,y) \memof
(\pi'_1,\pi'_2,\piev')$ it suffices to show $y=\piev' t$.  Now $\piev'
t = \piev \kappa t$. We have $\pi_1 \kappa t = \chi \pi'_1 t = \chi v
= \pi_1 s $ and $\pi_2 \kappa t = \pi'_2 t = \pi_2 s$. By the pullback
property $\pi_1$ and $\pi_2$ are jointly mono, so $\kappa t =s$. Hence
$y= \piev s = \piev \kappa t = \piev' t$ as desired. $\qed$

\begin{thm}  \label{th74}
Let $\cal C$ be a cartesian category satisfying (G). Let $Y \to^g X \to^f I$ be fixed morphisms.
Suppose that the  $\Pi$-diagram
\begin{equation} \label{diag74}
\bfig
\qtriangle(0,0)/<-`->`/[Y``;\piev`g`]         
\square(550,0)[P`F`X`I;\pi_1`\pi_2`\vf`f]
\efig
\end{equation}
is universal for $Y \to^g X \to^f I$. Then for every $i \in I$ and for every pair of morphisms 
$\psi=(r_1,r_2):R \to (X,Y)$ satisfying (A1) -- (A3), there is a unique $v \in F$  with $\vf v =i$ such that
for all $x \in X$ and $y \in Y$ 
\begin{equation} \label{equ74}
(x,y) \memof \psi  \Longleftrightarrow (v, x, y) \memof \alpha.
\end{equation}
Here $\alpha=(\pi_1, \pi_2, \piev): P \to (F,X,Y)$.
\end{thm}
{\flushleft \bf Proof.}  By Lemma \ref{lem21}  (\ref{defpid01}) is a
$\Pi$-diagram. Since (\ref{diag74}) is a universal diagram, there is a 
map $v:1 \to F$ such that $\vf v = i$ and  for all $u \in 1$, $x \in
X$ and $y\in Y$,
$$(u,x,y) \memof (!,r_1,r_2) \Longleftrightarrow (vu,x,y) \memof \alpha$$
(by Lemma \ref{lem22}). But $vu=v$ and $(u,x,y) \memof (!,r_1,r_2)$ is
equivalent to $(x,y) \memof \psi$, so (\ref{equ74}) is proved. $\qed$

\medskip
There is a converse

\begin{thm} \label{th75}
Let $\cal C$ be a cartesian category satisfying (G). Let $Y \to^g X \to^f I$ be fixed morphisms.
Consider the  $\Pi$-diagram
\begin{equation} \label{pidiag1}
\bfig
\qtriangle(0,0)/<-`->`/[Y``;\piev`g`]         
\square(550,0)[P`F`X`I;\pi_1`\pi_2`\vf`f]
\efig
\end{equation}
and let $\alpha=(\pi_1, \pi_2, \piev): P \to (F,X,Y)$.

Suppose that for every $i \in I$ and for every pair of morphisms 
$\psi=(r_1,r_2):R \to (X,Y)$ satisfying (A1) -- (A3), there is a unique $v \in F$  with $\vf v =i$ such that
for all $x \in X$ and $y \in Y$ 
$$(x,y) \memof \psi  \Longleftrightarrow (v, x, y) \memof \alpha.$$
Then (\ref{pidiag1}) is universal for $Y \to^g X \to^f I$.
\end{thm}
{\flushleft \bf Proof.} Let   
\begin{equation} 
\bfig                                                                          
\qtriangle/<-`->`/[Y``;\piev'`g`]         
\square(550,0)[P'`F'`X`I;\pi'_1`\pi'_2`\vf'`f]
\efig
\end{equation} 
be an arbitrary $\Pi$-diagram. For $v' \in F'$ form the pullback
\begin{equation}   \label{qpb1}
\bfig                                                                          
\square(0,0)[Q`1`P`F';!`q`v'`\pi'_1]
\efig
\end{equation} 
Then the composed diagram
\begin{equation}   \label{qpb2}
\bfig                                                                          
\qtriangle/<-`->`/[Y``;\piev'q`g`]         
\square(550,0)[Q`1`X`I;!`\pi'_2q`\vf' v'`f]
\efig
\end{equation} 
is, obviously, again a $\Pi$-diagram.  For $x \in X$ and $y \in Y$ we then have
$$(x,y) \memof (\pi'_2 q, \piev'q) \Longleftrightarrow (v',x,y) \memof
(\pi'_1,\pi'_2,\piev').$$ 
Indeed, suppose $x= \pi'_2 q u$ and $y=
\piev' q u$ for some $u\in Q$.  We have by (\ref{qpb1}) $v'= v' ! u =
\pi'_1 q u$. Hence $(v',x,y) \memof
(\pi'_1,\pi'_2,\piev')$. Conversely, suppose $v'= \pi'_1t$, $x= \pi'_2
t$ and $y= \piev' t$ for some $t \in P'$. From $v'= \pi'_1t$ it
follows by (\ref{qpb1}) that there is a unique $s \in Q$ with $t= q
s$. Thus $x= \pi'_2 q s$ and $y= \piev' q s$ and hence $(x,y) \memof
(\pi'_2 q, \piev'q)$.

Now (\ref{qpb2}) is a $\Pi$-diagram so $\psi = (\pi'_2 q, \piev'q): Q
\to (X,Y)$ satisfies (A1) -- (A3) for $i= \vf' v'$ (by Lemma
\ref{lem22}). Hence by assumption we have that there is a unique $v
\in F$ with $\vf v= \vf' v'$ and for all $x \in X$ and $y \in Y$
$$(v,x,y) \memof (\pi_1, \pi_2, \piev) \Longleftrightarrow (x,y) \memof \psi.$$

In conclusion, for every $v' \in F'$ there is a unique $v \in F$ such that $\vf v= \vf' v'$ and
\begin{equation} \label{eqv75}
(\forall x \in X)(\forall y \in Y)[(v,x,y) \memof (\pi_1, \pi_2,
    \piev) \Longleftrightarrow (v',x,y) \memof (\pi'_1, \pi'_2,
    \piev')].
\end{equation}

By unique choice (Proposition \ref{fungraph}) and Theorem \ref{elemcomp} there 
is $\chi:F' \to F$ so that for all $v'\in F$ it holds that
$\vf \chi v'= \vf' v'$ and
\begin{equation} \label{star35}
(\forall x \in X)(\forall y \in Y)[(\chi v',x,y) \memof (\pi_1, \pi_2, \piev) \Longleftrightarrow (v',x,y) \memof (\pi'_1, \pi'_2, \piev')].
\end{equation}
Hence according to Lemma \ref{lem22} the unique map $\kappa:P' \to P$ satisfying
$\pi_1 \kappa = \chi \pi'_1$ and  $\pi_2 \kappa = \pi'_2$ also satisfies $\piev \kappa = \piev'$.
To finish the proof we have to show that $\chi$  is unique. Suppose
that $\theta: F' \to F$ satisfies $\vf \theta = \vf'$ and that $\lambda:P' \to P$
is the unique map with $\pi_1 \lambda = \theta \pi'_1$ and $\pi_2
\lambda = \pi'_2$, and that this $\lambda$ satisfies $\piev \lambda = \piev'$. Then
by Lemma \ref{lem22}, it holds for all $v' \in F$
$$(\forall x \in X)(\forall y \in Y)[(\theta v',x,y) \memof (\pi_1, \pi_2, \piev) \Longleftrightarrow (v',x,y) \memof (\pi'_1, \pi'_2, \piev')].$$
Thus by the uniqueness in (\ref{eqv75})  for all $v' \in F'$, 
$\theta v' = \chi v'$. Hence $\theta = \chi$. $\qed$

\begin{cor}{\em (Theorem \ref{lccchar})}
Let $\cal C$ be a cartesian category satisfying (G). Then
$\cal C$ is locally cartesian closed if and only if $\cal C$ satisfies
the axiom ($\Pi$).
\end{cor}
{\flushleft \bf Proof.} ($\Rightarrow$) This is Theorem \ref{th74}. 
($\Leftarrow$) Suppose axiom ($\Pi$) holds. By Theorem \ref{th75} this
says that every $Y \to^g X \to^f I$ has a universal $\Pi$-diagram. 
$\qed$


{\sc
Erik Palmgren

Swedish Collegium for Advance Study, Uppsala, and

Institut Mittag-Leffler, Djursholm, and

Department of Mathematics, Uppsala University, Sweden.
}
\end{document}